\documentclass[a4paper,12pt]{article} 
\usepackage[english]{babel}
\usepackage{graphicx,latexsym,amsmath,amsthm,amssymb,color,indentfirst,xfrac}
 \usepackage[usenames,dvipsnames]{pstricks}
\usepackage{mathrsfs} 
 \usepackage{epsfig}
\usepackage{pst-grad} 
 \usepackage{pst-plot} 

\newtheorem{satz}{thm1}

\newtheorem{sat}{thm2}

\newtheorem{thm}[sat]{Theorem}
\newtheorem{lem}[satz]{Lemma}
\newtheorem{cor}[satz]{Corollary}
\newtheorem{prop}[satz]{Proposition}
\newtheorem{ex}[satz]{Example}
\theoremstyle{definition}
\newtheorem{defn}[satz]{Definition}

\newcommand{\be}{\begin{equation}}           
\newcommand{\ee}{\end{equation}}

\newcommand{\ba}{\begin{align}}                
\newcommand{\ea}{\end{align}}

\newcommand{\bal}{\begin{align*}}              
\newcommand{\eal}{\end{align*}}

\newcommand{\bxx}{\begin{ex}}

\newcommand{\exx}{\end{ex}}

\newcommand{\txsit}

\newenvironment{pr}
{\begin{trivlist}
\item[\hskip\labelsep{\bf Proof.}]}                     
{$\hfill\Box$\end{trivlist}}

\title{Sufficient conditions on planar graphs to have a relaxed DP-$3$-colorability}
\author {Pongpat Sittitrai\\ 
{\small\em Department of Mathematics, Faculty of Science, Khon Kaen University, 40002, Thailand }\\  
{\small\em E-mail address: pongpat\_s@kkumail.com} 
\and Kittikorn Nakprasit \footnote{Corresponding Author} \\ 
{\small\em Department of Mathematics, Faculty of Science, Khon Kaen University, 40002, Thailand }\\
{\small\em E-mail address: kitnak@hotmail.com}}
\date{}

\begin{document}

\maketitle

\begin{center}{\bf Abstract}\end{center}
\indent\indent

It is known that DP-coloring is a generalization of a list coloring in simple graphs 
and many results in list coloring can be generalized in those of DP-coloring. 
In this work, we introduce a relaxed DP-coloring which is a generalization 
if a relaxed list coloring. 
We also shows that every planar graph $G$ without $4$-cycles or $6$-cycles is DP-$(k,d)^*$-colorable. 
It follows immediately that $G$ is $(k,d)^*$-choosable.

\section{Introduction}
Every graph in this paper is finite, simple, and undirected. 
Embedding a graph $G$ in the plane, we let $V(G), E(G),$ and $F(G)$ 
denote the vertex set, edge set, and face set of $G.$ 
For $U \subseteq V(G),$ we let $G[U]$ denote 
the subgraph of $G$ induced by $U.$ For $X, Y \subseteq V(G)$ 
where $X$ and $Y$ are disjoint, we let $E_G(X,Y)$ be the set 
of all edges in $G$ with one endpoint in $X$ and the other in $Y.$ 

The concept of choosability was introduced by Vizing in 1976 \cite{Vizing} 
and by Erd\H os, Rubin, and Taylor in 1979 \cite{Erdos}, independently. 
A $k$-\emph{list assignment} $L$ of a graph $G$ assigns a list $L(v)$ 
(a set of colors) with $|L(v)|= k$ to each vertex $v.$ 
A graph $G$ is $L$-colorable if there is a proper coloring $f$ where $f(v)\in L(v).$   
If $G$ is $L$-colorable for every $k$-assignment $L,$ then we say  $G$ 
is $k$-\emph{choosable}. 

Dvo\v{r}\'{a}k and Postle \cite{DP} introduced a generalization 
of list coloring in which they called a \emph{correspondence coloring}. 
But following Bernshteyn, Kostochka, and Pron \cite{BKP},  
we call it a \emph{DP-coloring}.   


\begin{defn}\label{cover}
Let $L$ be an assignment of a graph $G.$ 
We call $H$ a \emph{cover} of $G$ 
if it satisfies all the followings:\\ 
(i) The vertex set of $H$ is $\bigcup_{u \in V(G)}(\{u\}\times L(u)) =
\{(u,c): u \in V(G), c \in L(u) \};$\\
(ii) $H[u\times L(u)]$ is a complete graph for every $u \in V(G);$\\
(iii) For each $uv \in E(G),$ 
the set $E_H(\{u\}\times L(u), \{v\}\times L(v))$ is a matching (maybe empty). 
(iv) If $uv \notin E(G),$ then no edges of $H$ connect  
$\{u\}\times L(u)$ and  $\{v\}\times L(v).$   
\end{defn}

\begin{defn}\label{DP} 
An $(H,L)$-coloring of $G$ is an independent set in 
a cover $H$ of $G$ with size $|V(G)|.$ 
We say that a graph is \emph{DP-$k$-colorable} if $G$ has 
an $(H,L)$-coloring for every $k$-assignment $L$ and every cover $H$ of $(G.$   
The \emph{DP-chromatic number} of $G,$ denoted by $\chi_{DP}(G),$ 
is the minimum number $k$ such that $G$ is DP-$k$-colorable. 
\end{defn}

If we define edges on $H$ to match exactly the same colors in $L(u)$ and $L(v)$ 
for each $uv \in E(H),$
then $G$ has an $(H,L)$-coloring if and only 
if $G$ is $L$-colorable.  
Thus DP-coloring is a generalization of list coloring. 
This also implies that  $\chi_{DP}(G) \geq \chi_l(G).$ 
In fact, the difference of these two chromatic numbers can be arbitrarily large. 
For graphs with average degree $d,$ Bernshteyn \cite{Bern} 
showed that   $\chi_{DP}(G) = \Omega(d /\log d),$   
while Alon \cite{Alon} showed that   $\chi_l(G) = \Omega(\log d).$ 

Dvo\v{r}\'{a}k  and Postle \cite{DP} showed that 
$\chi_{DP}(G) \leq 5$ for every planar graph $G.$  
This extends a seminal result by Thomassen \cite{Tho} on list colorings. 
On the other hand, Voigt \cite{Vo1}  gave an example of a planar graph 
which is not $4$-choosable (thus not DP-$4$-colorable). 
It is of interest to obtain sufficient conditions 
for planar graphs to be DP-$4$-colorable. 
Kim and Ozeki \cite{KimO} showed that planar graphs without $k$-cycles 
are DP-$4$-colorable for each $k =3,4,5,6.$  
Kim and Yu  \cite{KimY} extended the result on $3$- and $4$-cycles 
by showing that planar graphs without triangles adjacent to $4$-cycles  
are DP-$4$-colorable. 

\indent The concept of improper choosability was 
independently introduced by \v{S}krekovski \cite{skre}, 
and Eaton and Hull \cite{Eaton}. 
A graph $G$ is $(L,d)^*$-\emph{colorable} if there is a coloring $f$ where $f(v)\in L(v)$ 
such that every subgraph induced by vertices with the same color has maximum degree at most $d.$    
If $G$ is $(L,d)^*$-colorable for every $k$-assignment $L,$ then we say  $G$ 
is $(k,d)^*$-\emph{choosable}. 

In 1986, Cowen, Cowen, and Woodall \cite{ex4} constructed 
a planar graph that is not $(3,1)^*$-choosable. 
Many sufficient conditions for planar graphs 
to be $(3,1)^*$-choosable are studied. 
Zhang \cite{5and6} showed that every planar graphs 
without $5$-cycles or $6$-cycles is $(3,1)^*$-choosable. 
Chen and Raspaud \cite{P4adj4and3} proved that every planar graph 
without $4$-cycles adjacent to $3$- or $4$-cycles is $(3, 1)$*-choosable.  
Chen, Raspaud, and Wang \cite{Padjtriand6} proved that 
every planar graph without adjacent triangles or $6$-cycles 
is $(3, 1)$*-choosable. 

Inspired by $DP$-coloring, we define a generalization of a relaxed list coloring as follows.  

\begin{defn}\label{DPrelax} 
Let $H$ be a cover of a graph $G$ with respect to a list assignment $L.$ 
A $d$-\emph{representative} set of $G$ is a set of vertices $S$ 
in $H$ such that \\
(1) $|S|=|V(G)|,$ \\
(2) $u \neq v$ for any two different members $(u, c)$ and $(v,c')$ in $S,$ and\\   
(3) $H[S]$ has maximum degree at most $d.$ \\
An $(H,L,d)$-coloring of $G$ is a $d$-representative set.
We say that a graph is \emph{DP-$(k,d)^*$-colorable} if $G$ has an $(H,L,d)$-coloring for every $k$-assignment $L$ and every cover $H$ of $G.$  
Since only $d=1$ is considered in this paper, 
we write a representative set instead of a $1$-representative set. 
\end{defn}

If we define edges on $H$ to match exactly the same colors in $L(u)$ and $L(v)$ 
for each $uv \in E(H),$
then $G$ has an $(H,L,d)$-coloring if and only 
if $G$ is $(L,d)^*$-colorable. 
This follows immediately that $G$ is DP-$(k,d)^*$-colorable implies $G$ is  $(k,d)^*$-choosable.   

In this work, we have the following result. 

\begin{thm}\label{main2}
Every planar graph without $4$-cycles or $6$-cycles is DP-$(3,1)^*$-colorable.
\end{thm}


\section{Structure Obtained from Condition on Cycles}
First, we introduce some notations and definitions. 
A $k$-vertex ($k^+$-vertex, $k^-$-vertex, respectively) is 
a vertex of degree $k$ (at least $k,$ at most $k,$ respectively). 
The same notations are applied to faces. 
A \emph{$(d_1,d_2,\dots,d_k)$-face} $f$ is a face of degree $k$ 
where  vertices on $f$ have degree $d_1,d_2,\dots,d_k$ in  a cyclic order. 
A \emph{$(d_1,d_2,\dots,d_k)$-vertex} $v$ is a vertex of degree $k$ 
where  faces incident to $v$ have degree $d_1,d_2,\dots,d_k$ in a cyclic  order.  
A $(3,4^+,4^+)$-face $f$ is called a \emph{pendant} $3$-face of $v$ if $v$ is not a vertex of $f$ but adjacent to a $3$-vertex of $f.$

Let $G$ be a graph without $4$-cycles or  $6$-cycles. 
The following property is straightforward. 
\begin{prop}\label{prop1} 
	A $3$-face in $G$ does not share exactly one edge with $6^-$-faces. 
\end{prop}
Proposition \ref{prop1} yields the following two Propositions. 
\begin{prop}\label{prop2} 
If $f$ is a pendant $3$-face of $v$ where $u$ is a $3$-vertex in $f$, 
then  two  faces that are  incident to both $v$ and $u$ are  $7^+$-faces. 
\end{prop}
\begin{prop} \label{prop11} 
 Every vertex $v$ is incident to at most $\lfloor\frac{d(v)}{2}\rfloor$ $3$-faces.
\end{prop} 

\section {Structure of Minimal Non DP-(3,1)-colorable Graphs}
\begin{defn}\label{residual} 
	Let $H$ be a cover of $G$ with a list assignment $L.$ 
	Let $G'= G-F$ where $F$ is an induced subgraph of $G.$ 
	A list assignment $L'$ is a \emph{restriction of $L$} on $G'$ 
	if $L'(u) = L(u)$ for each vertex in $G'.$ 
	A graph $H'$ is a \emph{restriction of $H$} on $G'$ 
	if $H'= H[\{\{v\} \times L(v): v \in V(G')\}.$ 
	Assume $G'$ has an $(H',L',1)$-coloring with a representative set $R'$ in $H'$ such that  
	$|R'|= |V(G)|-|V(F)|.$\\
	\indent A \emph{residual list assignment} $L^*$ of $F$ is defined by    
	$$L^*(x)=L(x)-\bigcup_{ux\in E(G)}\{c'\in L(x) : 
	(u,c)(x,c')\in E(H) \text{ and }(u,c)\in I'\}$$
	for each $x \in V(F).$ \\
	\indent A \emph{residual cover $H^*$} is defined by $H^*=H[\{\{x\} \times L^*(x): x \in V(F)\}].$
\end{defn}
From above definitions, we have the following fact. 

\begin{lem}\label{extend} 
	Assume $G$ has an induced subgraph $G'$ and a cover $H$ with a list assignment $L.$ 
	Let $H'$ be a restriction of $H$ on $G'.$ 
    If $G'$ has an $(H',L',1)$-coloring with a representative set $R'$ in $H'$ such that  
	$|R'|= |V(G)|-|V(F)|,$ then a residual cover $H^*$ is a cover of $F$ with an assignment $L^*. $
	Furthermore, if $F$ is $(H^*,L^*,1)$-colorable, then $G$ is $(H,L,1)$-colorable. 	
\end{lem}
\begin{pr}
	One can check from the definitions of a cover and a residual cover that 
	$H^*$ is a cover of $F$ with an assignment $L^*.$\\
	\indent Suppose that  $F$ is $(H^*,L^*,1)$-colorable. 
	Then $H^*$ has a representative set $R^*$ with $|R^*|= |F|.$ 
	It follows from Definition \ref{residual} that no edges connect $H^*$ and $R'.$ 
	Additionally, $R'$ and $R^*$ are disjoint. 
	Altogether, we have that $R=R'\cup R^*$ is a representative set in $H$   
	with $|R|=(|V(G)|-|V(F))+|V(F)|=|V(G)|.$ 
	Thus $G$ is $(H,L,1)$-colorable. 	
\end{pr}
From now on, let  $G$ be a minimal non DP-$4$-colorable graph. 

\begin{lem}\label{lem1} 
	Each vertex in $G$ is a $3^+$-vertex.
\end{lem}
\begin{pr} 
	Suppose to the contrary that $G$ has a vertex $x$ degree at most $2.$ 
	Let $L$ be a $3$-assignment 
	and let $H$ be a cover of $G$ such that $G$ has no $(H,L,1)$-coloring.       
	By the minimality of $G,$ the subgraph $G'=G-x$ admits $(H',L',1)$-coloring 
	where $L'$ (and $H'$) is a restriction of $L$ (and $H$, respectively) in $G'.$ 
	Thus there is a representative set $R'$ with $|R'|= |G'|$ in $H'.$  
	Consider a residual list assignment $L^*$ on $x.$ 
	Since $|L(x)| = 3$ and $d(x)\leq 2$, we obtain $|L^*(x)|\geq 1.$  
	Clearly, $\{(x,c)\}$ where $c \in L^*(x)$ is a representative set in $G[\{x\}].$  
	Thus $G[\{x\}]$ is $(H^*,L^*,1)$-colorable. 
	It follows from Lemma \ref{extend} that the graph $G$ is $(H,L,1)$-colorable, a contradiction. 
\end{pr}

\begin{lem}\label{lem2} 
	Each neighbor of $3$-vertex in $G$ is a $4^+$-vertex.
\end{lem}

\begin{pr} Suppose to the contrary that there are adjacent $3$-vertices $u$ and $v.$  
	Let $L$ be a $3$-assignment 
	and let $H$ be a cover of $G$ such that $G$ has no $(H,L,1)$-coloring.   
	By the minimality of $G,$ 
	the subgraph $G'=G-\{u,v\}$ admits an $(H',L',1)$-coloring 
	where $L'$ (and $H'$, respectively) is a restriction of $L$ (and $H$, respectively) in $G'.$ 
	Thus there is a representative set $R'$ with $|R'|= |G'|$ in $H'.$  
	Consider a residual list assignment $L^*$ on $G[\{u,v\}].$ 
	We have $|L^*(u)|$ and  $|L^*(v)| \geq 1$. 
	Clearly, $\{(v,c),(u,c')\}$ where $c \in L^*(v)$ and  $c' \in L^*(u)$  is a representative set in $G[\{u,v\}].$ 
	We obtain a representative set $R^*$ with $|R^*|=2.$	
	Thus $G[\{u,v\}]$ is $(H^*,L^*,1)$-colorable. 
	It follows from Lemma \ref{extend} that  $G$ is $(H,L,1)$-colorable, a contradiction. 
\end{pr}
\begin{lem}\label{lem3} 
	Each  $4$-vertex in $G$ is adjacent to at most two   $3$-vertices.
\end{lem}

\begin{pr} Suppose to the contrary that a $3$-vertex  $v$ is adjacent to  three $3$-vertices, $u_1,u_2,$ and $u_3$.  
	It follows from Lemma \ref{lem2} that $u_i$ is not adjacent to $u_j$ for $i,j\in\{1,2,3\}.$ 
	Let $L$ be a $3$-assignment 
	and let $H$ be a cover of $G$ such that $G$ has no $(H,L,1)$-coloring.   
	By the minimality of $G,$ 
	the subgraph $G'=G-\{v,u_1,u_2,u_3\}$ admits an $(H',L',1)$-coloring 
	where $L'$ (and $H'$, respectively) is a restriction of $L$ (and $H$, respectively) in $G'.$ 
	Thus there is a representative set $R'$ with $|R'|= |G'|$ in $H'.$  
	Consider a residual list assignment $L^*$ on $G[\{v,u_1,u_2,u_3\}].$ 
	Since $|L(v)|= 3$ for every $v\in V(G),$ 
	we have $|L^*(v)|\geq2$ and  $|L^*(u_i)| \geq 1$ for $i\in\{1,2,3\}$. 
Let $H^*$ be an residual cover of $G[\{v,u_1,u_2,u_3\}].$ 
	First of all, we  choose a color $c_i$ from  $L^*(u_i)$. So there is a
color $c$ in $L^*(v)$ such that $(v,c)$ is adjacent to at most one of $\{(u_i,c_i)\}$ for $i\in\{1,2,3\}$. 
Clearly, $\{(v,c),(u_1,c_1),(u_2,c_2),(u_3,c_3)\}$ where $c \in L^*(v)$ and  $c_i \in L^*(u_i)$ for $i\in\{1,2,3\}$   is a representative set in $G[\{v,u_1,u_2,u_3\}].$  
Thus	we obtain a representative set $R^*$ with $|R^*|=4=|G[\{v,u_1,u_2,u_3\}]|.$ 
	Thus $G[\{v,u_1,u_2,u_3\}]$ is $(H^*,L^*,1)$-colorable. 
	It follows from Lemma \ref{extend} that  $G$ is $(H,L,1)$-colorable, a contradiction. 
\end{pr}
From Lemma \ref{lem2}, we obtain the upper bound of the number of incident $3$-vertices of a face in $G.$  
\begin{cor}\label{cor4}  
	Each faces $f$ in $G$ is incident to at least $\frac{d(f)}{2}$ $3$-vertices. 
\end{cor}
\section{Main Result}
\begin{thm} \label{thm1} Every planar graph without $4$-cycles or $6$-cycles is DP-$(3,1)^*$-colorable.
\end{thm}
\begin{pr}
\indent Suppose that $G$ is a minimal counterexample. Then each vertex in $G$ is a $3^+$-vertex by Lemma \ref{lem1}. The discharging process is as follows. Let the initial charge of a vertex $u$ in $G$ be $\mu(u)=2d(u)-6$ and the initial charge of a face $f$ in $G$ be $\mu(f)=d(f)-6$. Then by Euler's formula $|V(G)|-|E(G)|+|F(G)|=2$ and by the Handshaking lemma, we have
$$\displaystyle\sum_{u\in V(G)}\mu(u)+\displaystyle\sum_{f\in F(G)}\mu(f)=-12.$$
\indent Now, we establish a new charge $\mu^*(x)$ 
for all $x\in V(G)\cup F(G)$ by transferring charge from one element to another
and the summation of new charge $\mu^*(x)$ remains $-12$. 
If the final charge  $\mu^*(x)\geq 0$ for all $x\in V(G)\cup F(G)$, 
then we get a contradiction and the proof is completed.\\
\indent The discharging rules are\\
(R1) Every $4^+$-vertex sends charge $1$ to each incident $3$-face.\\
(R2) Every $4^+$-vertex sends charge $\frac{1}{3}$ to each incident $5$-face.\\
(R3) Every $4^+$-vertex sends charge $\frac{1}{3}$ to each pendent  $3$-face.\\
(R4) Every $7^+$-face sends charge $\frac{1}{3}$ to each incident $3$-vertex.\\
(R5) Every $3$-vertex sends charge $\frac{2}{3}$ to each incident  $3$-face.\\
\indent Next, we  show that the final charge  $\mu^*(u)$  is nonnegative.\\
\indent \textit{CASE 1:} A $3$-vertex $v.$\\
\indent If $v$ is not incident to any $3$-face, then $\mu^*(v) \geq  0.$ 
 If $v$ is incident to a $3$-face, then $v$ is a $(3,7^+,7^+)$-vertex by Proposition \ref{prop1}. Thus 
 $\mu^*(v) \geq \mu(v) -\frac{2}{3}+2\times\frac{1}{3}= 0$ by (R4) and (R5).\\
\indent \textit{CASE 2:} A $4$-vertex $v.$\\
\indent It follows from by Proposition \ref{prop11} that  $v$ is incident to at most two $3$-faces.  
If $v$ is incident to two $3$-faces,  
then$v$ is a $(3,7^+,3,7^+)$-vertex by Proposition \ref{prop1} and $v$ has no any pendent $3$-faces by  Lemma \ref{prop2}. 
Thus $\mu^*(v) \geq \mu(v) -2\times1= 0$ by (R1).
If $v$ is incident to exactly one $3$-face, then $v$ is a $(3,7^+,5^+,7^+)$-vertex by Proposition \ref{prop1}. 
Moreover, $v$ has at most two pendent $3$-faces by  Lemma \ref{lem3}.
Thus $\mu^*(v) \geq \mu(v) -1-3\times\frac{1}{3}= 0$ by (R1), (R2), and (R3).
If $v$ is not incident to any $3$-face, then $v$ has at most two pendent  $3$-faces by Lemma \ref{lem3}. 
Thus $\mu^*(v) \geq  \mu(v) -6\times\frac{1}{3}= 0$ by (R2) and (R3).\\
\indent \textit{CASE 3:} A $5^+$-vertex $v.$\\
\indent 
To facilitate the calculation, we redefine the discharging rule for $v$ and its incident faces $f_1,f_2,\dots,f_{d(v)}.$ 
Let $v$ send charge $\frac{2}{3}$ to each incident face. 
We have  $\mu^*(v) \geq (2d(v)-6) - d(v)\times\frac{2}{3} \geq 0.$    
Now, let each non $3$-face $f_i$ send charge $\frac{1}{6}$ to each adjacent $3$-face of $v$ and each adjacent pendent $3$-face of $v$. That means the remaining charge is received by $v$ as follows,\\ 
(1) Each $3$-face of $v$ receive charge at least  $\frac{2}{3}+2\times\frac{1}{6}=1$ by Proposition \ref{prop1}.\\ 
(2) Each  $5^+$-face of $v$ receive charge at least  $\frac{2}{3}-2\times\frac{1}{6}=\frac{1}{3}$ by Proposition \ref{prop2}.\\ 
(3) Each pendant $3$-face of $v$ receive charge at least  $2\times\frac{1}{6}=\frac{1}{3}$ by Proposition  \ref{prop2}.\\
 \indent One can see that charge of each $f_i$ is at least that obtains from (R1), (R2), and (R3). 
 Thus $\mu^*(v) \geq 0.$\\
\indent \textit{CASE 4:} A $3$-face $f.$\\
\indent It follows from  Lemma \ref{lem2} that $f$ is a $(3^+,4^+,4^+)$-face. 
If $f$ is a $(3,4^+,4^+)$-face, then  $f$ is a pendant $3$-face of some $4^+$-vertex by Lemma \ref{lem2}. 
Thus $\mu^*(f) \geq  \mu(f) +2\times1+\frac{2}{3}+\frac{1}{3}= 0$ by (R1), (R3), and (R5). 
If $f$ is a $(4^+,4^+,4^+)$-face, then  $\mu^*(f) \geq  \mu(f) +3\times1= 0$ by (R1).\\
  \indent \textit{CASE 5:} A $5^+$-face $f.$\\
  \indent If $f$ is a $5$-face, then $f$ is incident to at least three  $4^+$-vertices by Corollary \ref{cor4}.  
 Thus  $\mu^*(f) \geq  \mu(f) +3\times\frac{1}{3}= 0$ by (R2).
If $f$ is a $7$-face, then $f$ is incident to at most three $3$-vertices by Corollary \ref{cor4}. 
Thus $\mu^*(f) \geq  \mu(f) -3\times\frac{1}{3}= 0$ by (R4). 
If $f$ is a $8^+$-face, then $f$ is incident to at most $\frac{d(f)}{2}$ $3$-vertices by Corollary \ref{cor4}. 
Thus $\mu^*(f) \geq  \mu(f) -\frac{d(f)}{2}\times\frac{1}{3}> 0$ from (R4).\\
\indent This completes the proof. 
\end{pr}
\section{Acknowledgments}
\indent The first author is supported by Development and Promotion of Science and Technology
talents project (DPST).

\end{document}